\overfullrule=0pt
\centerline {\bf A property of strictly convex functions}
\centerline  {\bf which differ from each other by a constant on the boundary of their domain}\par
\bigskip
\bigskip
\centerline {BIAGIO RICCERI}\par
\bigskip
\bigskip
 In this paper, we prove, in particular, the following result:  Let $E$ be a reflexive real Banach space and let $C\subset E$ be a closed convex set, with non-empty interior, whose boundary is sequentially weakly closed and non-convex. Then, for every function $\varphi:\partial C\to {\bf R}$ and for every convex set $S\subseteq E^*$ dense in $E^*$, there exists $\tilde\gamma\in S$ having the following property: for
every strictly convex lower semicontinuous function $J:C\to {\bf R}$, G\^ateaux differentiable in $\hbox {\rm int}(C)$, such that $J_{|\partial C}-\varphi$ is constant
in $\partial C$ and $\lim_{\|x\|\to +\infty}{{J(x)}\over {\|x\|}}=+\infty$ if $C$ is unbounded,
 $\tilde\gamma$ is an algebraically
interior point of $J'(\hbox {\rm int}(C))$ (with respect to $E^*$).\par
\par
\bigskip
\bigskip
{\bf 1. Introduction}\par
\bigskip
When $E$ is a topological space, a function $f:E\to {\bf R}$ is said to be inf-compact if, for every $r\in {\bf R}$, the set $f^{-1}(]-\infty,r])$ is compact. When $E$ is a real
vector space and $A\subseteq E$, a point $x_0\in A$ is said to be an algebraically interior point of $A$ (with respect to $E$) if, for every $y\in E$, there exists $\delta>0$ such that
$x_0+\lambda y\in A$ for all $\lambda\in [0,\delta]$. The algebraic interior of $A$ is the set of all its algebraically interior points.\par
\smallskip
After recalling these two notions, we can state the result which the title refers to:\par
\medskip
THEOREM 1.1. - {\it Let $E$ be a reflexive real Banach space and let $C\subset E$ be a closed convex set, with non-empty interior, whose boundary is sequentially weakly closed and non-convex.\par
Then, for every function $\varphi:\partial C\to {\bf R}$ and for every convex set $S\subseteq E^*$ dense in $E^*$, there exists $\tilde\gamma\in S$ having the following property: for
every strictly convex lower semicontinuous function $J:C\to {\bf R}$, G\^ateaux differentiable in $\hbox {\rm int}(C)$, such that $J_{|\partial C}-\varphi$ is constant
in $\partial C$ and $\lim_{\|x\|\to +\infty}{{J(x)}\over {\|x\|}}=+\infty$ if $C$ is unbounded,
 $\tilde\gamma$ is an algebraically
interior point of $J'(\hbox {\rm int}(C))$ (with respect to $E^*$)}.\par
\medskip
Now, some comments are in order. The main feature of Theorem 1.1 is the fact that $\tilde\gamma$ does not depend on $J$. But, for a moment, consider simply the following by-product of Theorem 1.1:\par
 \smallskip
{\it If $C$ is as before, then, for every stricly convex lower semicontinuous function
$J:C\to {\bf R}$, G\^ateaux differentiable in $\hbox {\rm int}(C)$ and with $\lim_{\|x\|\to +\infty}{{J(x)}\over {\|x\|}}=+\infty$ if $C$ is unbounded
, the algebraic interior of $J'(\hbox {\rm int}(C))$ (with respect to $E^*$) is non-empty}.\par
\smallskip
 As far as we know, such a 
corollary itself is new when $E$ is infinite-dimensional. To the contrary, if $E$ is finite-dimensional, it can essentially be considered known, at least
if $J'$ is continuous. Indeed, in that case, $J'$, being injective (due to the strict convexity of $J$), 
turns out to be open, thanks to the invariance of domain theorem. The assumptions that $\partial C$ is sequentially weakly closed and that $\lim_{\|x\|\to +\infty}{{J(x)}\over {\|x\|}}=+\infty$ if $C$ is unbounded cannot be removed. Actually, consider the following
situation. Let $f:{\bf R}\to {\bf R}$ be a continuous, increasing and bounded function.
Define the functional $J:L^2([0,1])\to {\bf R}$ by
$$J(u)=\int_0^1\left (\int_0^{u(x)}f(t)dt\right)dx$$
for all $u\in L^2([0,1])$. Clearly, $J$ is stricly convex and $C^1$, with 
$$J'(u)=f\circ u$$
for all $u\in L^2([0,1])$ (after identifiyng $(L^2([0,1])^*$ to $L^2([0,1])$).
Notice that, since $J'(L^2([0,1]))\subseteq L^{\infty}([0,1])$, for each $A\subseteq L^2([0,1])$, the algebraic interior
of $J'(A)$ (with respect to $L^2([0,1]))$ is empty. Now, let $C$ be any closed ball in $L^2([0,1])$. So, the restriction of $J$ to $C$ is weakly inf-compact. In this case, the conclusion of the corollary fails since $\partial C$ is not sequentially weakly closed. On the other hand, when $C$ is as in the corollary, the conclusion fails since the restriction of $J$ to $C$ is not weakly inf-compact.\par
\smallskip
Now, come back to the full statement of Theorem 1.1. We observe that the assumption about the non-convexity of $\partial C$ cannot be removed. In this connection, let
$U$ be the set of all $C^1$ strictly convex functions $J:[0,+\infty[\to {\bf R}$, with $J(0)=0$ and $\lim_{x\to +\infty}{{J(x)}\over
{x}}=+\infty$. Fix $y\in {\bf R}$. Let $\mu>\max\{0,y\}$. Then, taking $J(x)=\mu(e^x-1)$, we have 
 $J\in U$ and $J'(x)>y$ for all $x>0$. So, taking $C=[0,+\infty[$ and $\varphi(0)=0$, the conclusion of the Theorem 1.1 does not hold:
the boundary of $C$ is convex.\par
\smallskip
Our proof of Theorem 1.1 is fully based on a minimax result obtained in [2]. More precisely, in the next section, we will establish a general
abstract result (just based on [2]) from which we will derive a series of consequences among which there is Theorem 1.1.\par
\bigskip
{\bf 2. Results}\par
\bigskip
First, let us recall the following result from [2]:\par
\medskip
THEOREM 2.A [2], Theorem 1.1). - {\it Let $X$ be a topological space, $F$ a topological vector
space, $Y\subseteq F$ a non-empty convex set 
and $g:X\times Y\to
{\bf R}$ a function satisfying the following conditions:\par
\noindent
$(a)$\hskip 5pt for each $y\in Y$, the function $g(\cdot,y)$ is lower semicontinuous and
inf-compact;\par
\noindent
$(b)$\hskip 5pt for each $x\in X$, the function $g(x,\cdot)$ is continuous and quasi-concave.\par
\noindent
$(c)$\hskip 5pt $\sup_Y\inf_Xg<\inf_X\sup_Yg$\ .\par
Then, there exists $y^*\in Y$ such that the function $g(\cdot,y^*)$ has at least two global minima.}\par
\medskip
In what follows, $E$ is a topological space and  $Y$ is a convex set in a topological vector space.\par
\smallskip
Let us introduce the two main classes of functions we deal with.
\smallskip
Let $X\subseteq C\subseteq E$ and let $\varphi:X\to {\bf R}$ be a given function. We denote by
${\cal B}(X,C,\varphi)$ the class of all functions $J:C\to {\bf R}$ such that $J_{|X}-\varphi$ is
constant in $X$.\par
\smallskip
Let $C\subseteq E$ and $S\subseteq Y$. Let $f:C\times S\to {\bf R}$ be a given function. We denote by
${\cal C}(C,S,f)$ the class of all functions $J:C\to {\bf R}$ such that, for each $y\in S$, the function
$f(\cdot,y)+J(\cdot)$ has at most one global minimum in $C$.\par
\smallskip
Our main abstract result is as follows:\par
\medskip
THEOREM 2.1. - {\it Let $X\subseteq C\subseteq E$, let $S\subseteq Y$ be a convex set dense in $Y$, let $f:C\times Y\to {\bf R}$ and let
$\varphi:X\to {\bf R}$. Assume that:\par
\noindent
$(a)$\hskip 5pt for each $y\in S$, the function $f(\cdot,y)+\varphi(\cdot)$ is lower semicontinuous and 
inf-compact in $X$\ ;\par
\noindent
$(b)$\hskip 5pt for each $x\in X$, the function $f(x,\cdot)$ is quasi-concave and continuous in $Y$\ ;\par
\noindent
$(c)$\hskip 5pt $\sup_Y\inf_X(f+\varphi)<\inf_X\sup_Y(f+\varphi)\ .$\par
Then, there exists a point  $y^*\in S$ such that 
$$\inf_{x\in C}(f(x,y^*)+J(x))<\inf_{x\in X}(f(x,y^*)+J(x))$$
for every $J\in {\cal B}(X,C,\varphi)\cap {\cal C}(C,S,f)$\ .}\par
\smallskip
PROOF. Consider the function $g:X\times S\to {\bf R}$ defined by
$$g(x,y)=f(x,y)+\varphi(x)$$
for all $(x,y)\in X\times S$. For each $x\in X$, by continuity of $f(x,\cdot)$ and density of $S$, we have
$$\sup_{y\in S}f(x,y)=\sup_{y\in Y}f(x,y)\ ,$$
and so
$$\sup_{y\in S}g(x,y)=\sup_{y\in Y}g(x,y)\ .$$
Therefore, in view of $(c)$, we have
$$\sup_{S}\inf_Xg\leq \sup_Y\inf_Xg<\inf_X\sup_Yg=\inf_X\sup_Sg\ .$$
The function $g$ is lower semicontinuous and inf-compact in $X$, and quasi-concave and continuous in $S$. So, thanks to Theorem 2.A,
there exists $y^*\in S$ such that the function
$(f(\cdot,y^*))_{|X}+\varphi(\cdot)$ has at least two global minima in $X$. Now, fix $J\in {\cal B}(X,C,\varphi)\cap {\cal C}(C,S,f)$.
Since $J_{|X}-\varphi$ is constant in $X$, the functions $(f(\cdot,y^*)+J(\cdot))_{|X}$ and $(f(\cdot,y^*))_{|X}+\varphi(\cdot)$ have the same
global minima in $X$. Arguing by contradiction, assume that
$$\inf_{x\in C}(f(x,y^*)+J(x))=\inf_{x\in X}(f(x,y^*)+J(x))\ .\eqno{(2.1)}$$
We know that $(f(\cdot,y^*)+J(\cdot))_{|X}$ has at least two global minima. But, in view of $(2.1)$, they turn out to be global minima
of $f(\cdot,y^*)+J(\cdot)$ in $C$, against the fact that $J\in {\cal C}(C,S,f)$. The proof is complete.\hfill $\bigtriangleup$\par
\medskip
Here is a remarkable corollary of Theorem 2.1.\par
\medskip
THEOREM 2.2. - {\it Let $X\subseteq C\subseteq E$, let $S\subseteq Y$ be a convex set dense in $Y$, let $f:C\times Y\to {\bf R}$ and let
$\varphi:X\to {\bf R}$. Assume that:\par
\noindent
$(i)$\hskip 5pt for each $y\in S$, the function $f(\cdot,y)+\varphi(\cdot)$ is lower semicontinuous and 
inf-compact in $X$\ ;\par
\noindent
$(ii)$\hskip 5pt for each $x\in X$, the function $f(x,\cdot)$ is quasi-concave and continuous in $Y$\ ;\par
\noindent
$(iii)$\hskip 5pt $\inf_X\sup_Yf=+\infty$ and there exists a finite set $A\subset X$ such that $\sup_Y\inf_Af<+\infty$\ .\par
Then, there exists a point  $y^*\in S$ such that 
$$\inf_{x\in C}(f(x,y^*)+J(x))<\inf_{x\in X}(f(x,y^*)+J(x))$$
for every $J\in {\cal B}(X,C,\varphi)\cap {\cal C}(C,S,f)$\ .}\par
\smallskip
PROOF. After observing that
$$\sup_Y\inf_X(f+\varphi)\leq \sup_Y\inf_Af+\sup_A\varphi<+\infty=\inf_X\sup_Yf=
\inf_X\sup_Y(f+\varphi)\ ,$$
the conclusion follows directly from Theorem 2.1.\hfill $\bigtriangleup$\par
\medskip
In the next two results, $E$ is also a real vector space (and the topology on $E$ is still arbitrary).\par
\medskip
THEOREM 2.3. - {\it Let $X\subseteq E$, let $F$ be a real normed space, let $S\subseteq F^*$ be a convex set weakly-star dense in $F^*$
let $I:{\rm conv}(X)\to {\bf R}$, let $\psi:{\rm conv}(X)\to F$ and let $\varphi:X\to {\bf R}$. Assume that that $\psi(X)$ is not convex and that
the function $(I+\eta\circ \psi)_{|X}+\varphi$ is lower semicontinuous and inf-compact in $X$ for
all $\eta\in S$.\par
Then, there exists $\tilde \eta\in S$ having the following property: for every function $J:{\rm conv}(X)\to {\bf R}$ such that
$J_{|X}-\varphi$ is constant in $X$ and $I+J+\eta\circ\psi$ is strictly convex for all $\eta\in S$, one has
$$\inf_{x\in \hbox {\rm conv}(X)}(I(x)+J(x)+\tilde\eta(\psi(x)))<\inf_{x\in X}(I(x)+J(x)+\tilde\eta(\psi(x)))\ .$$}\par
\smallskip
PROOF. Fix $y_0\in \hbox {\rm conv}(\psi(X))\setminus \psi(X)$. We apply Theorem 2.2 with $C=\hbox {\rm conv}(X)$ and $Y=F^*$, considering
$F^*$ equipped with the weak-star topology and taking
$$f(x,\eta)=I(x)+\eta(\psi(x)-y_0)$$
for all $(x,\eta)\in C\times Y$. Clearly, $f$ satisfies conditions $(i)$ and $(ii)$. Moreover, if $y_0=\sum_{i=1}^n\lambda_i\psi(x_i)$, where
$x_i\in X$, $\lambda_i\geq 0$, $\sum_{i=1}^n\lambda_i=1$, 
by Proposition 2.1 of [3], we know that
$$\sup_{\eta\in Y}\inf_{1\leq i\leq n}f(x_i,\eta)\leq \max_{1\leq i\leq n}I(x_i)<+\infty\ .$$
On the other hand, for each $x\in X$, since $\psi(x)\neq y_0$, we have
$$\sup_{\eta\in Y}\eta(\psi(x)-y_0)=+\infty\ ,$$
and so condition $(iii)$ is satisfied too. Now, Theorem 2.2 ensures the existence of $\tilde\eta\in S$ such that,
for every $J\in {\cal B}(X,C,\varphi)\cap {\cal C}(C,S,f)$, one has
$$\inf_{x\in \hbox {\rm conv}(X)}(f(x,\tilde\eta)+J(x))<\inf_{x\in X}(f(x,\tilde\eta)+J(x))$$
which means
$$\inf_{x\in \hbox {\rm conv}(X)}(I(x)+J(x)+\tilde\eta(\psi(x)))<\inf_{x\in X}(I(x)+J(x)+\tilde\eta(\psi(x)))\ .$$
To finish the proof, it is enough to remark that if $J:{\rm conv}(X)\to {\bf R}$ is such that $I+J+\eta\circ\psi$ is strictly convex for all $\eta\in S$, then $J\in {\cal C}(C,S,f)$.\hfill $\bigtriangleup$\par
\medskip
The following result is a particularly simple consequence of Theorem 2.3.\par
\medskip
THEOREM 2.4. - {\it  Let $X\subset E$ be a compact set, let $F$ be a real normed space,
and let $\psi:{\rm conv}(X)\to F$ be an affine operator, continuous
with respect to the weak topology on $F$, such that $\psi(X)$ is not convex, and let $\varphi:X\to {\bf R}$ be a lower semicontinous
function.\par
 Then, for every convex set $S\subseteq F^*$ weakly-star dense in $F^*$,
 there exists $\tilde\eta\in S$ having the following property: for every strictly convex function $J:\hbox {\rm conv}(X)\to {\bf R}$ such that
$J_{|X}-\varphi$ is constant in $X$, one has
$$\inf_{x\in \hbox {\rm conv}(X)}(J(x)+\tilde\eta(\psi(x)))<\inf_{x\in X}(J(x)+\tilde\eta(\psi(x)))\ .$$}\par
\smallskip
PROOF.  For each $\eta\in F^*$, the function
$\eta\circ\psi$ is continuous since $\eta$ is weakly continuous. So, $\eta\circ\psi+\varphi$ is lower semicontinuous and inf-compact, since
$X$ is compact. Hence, the assumptions of Theorem 2.3 are satisfied, with $I=0$. Now, our conclusion follows from that of Theorem 2.3 taken into account that if $\eta\in F^*$ and $J:\hbox {\rm conv}(X)\to {\bf R}$ is a strictly convex function, then so is $\eta\circ\psi+J$ since
$\psi$ is affine.\hfill $\bigtriangleup$\par
\medskip
The next consequence of Theorem 2.3 can be considered as the central result: actually, Theorem 1.1 is a corollary of it.\par
\medskip
THEOREM 2.5. - {\it Let $E$ be a reflexive real Banach space, let $C\subset E$ be a closed convex set, with non-empty interior, such that
$\partial C$ is sequentially weakly closed and non-convex, let $S\subseteq E^*$ be a convex set dense in $E^*$,
let $I:C\to {\bf R}$ be G\^ateaux differentiable in $\hbox {\rm int}(C)$, and let $\varphi:\partial C\to {\bf R}$.\par
Then, there exists $\tilde\gamma\in S$ having the following property: 
for every function $J:C\to {\bf R}$, G\^ateaux differentiable in $\hbox {\rm int}(C)$, such that
$J_{|\partial C}-\varphi$ is constant in $\partial C$ and $I+J$ is lower semicontinous and strictly convex,
with $\lim_{\|x\|\to +\infty}{{I(x)+J(x)}\over {\|x\|}}=+\infty$ if $C$ is unbounded,
and for every sequentially weakly lower semicontinuous function $G:C\to {\bf R}$, G\^ateaux differentiable in $\hbox {\rm int}(C)$,
there exists $\epsilon>0$ such that, for each $\lambda\in [0,\epsilon]$, the equation
$$I'(x)+J'(x)+\lambda G'(x)=\tilde\gamma$$
has at least one solution in $\hbox {\rm int}(C)$.}\par
\smallskip
PROOF. We apply Theorem 2.3 considering $E$ equipped with the weak topology (but the interior of $C$ is referred to the strong topology) and
taking $X=\partial C$, $F=E$ and $\psi(x)=x$ for all $x\in C$ (notice that $\hbox {\rm conv}(\partial X)\subseteq C$).
Of course, it is implicitly understood that there are functions $J_0:C\to {\bf R}$ such that
$J_{0_{|\partial C}}-\varphi$ is constant in $\partial C$, and $I+J_0$ is lower semicontinuous and strictly convex, with 
$\lim_{\|x\|\to +\infty}{{I(x)+J_0(x)}\over
{\|x\|}}=+\infty$ if $C$ is unbounded. So, if $J_0$ is such a function, it follows that,  for each $\eta\in E^*$, the function $I+J_0+\eta$ is weakly inf-compact in $C$ (in particular, notice that $\lim_{\|x\|\to +\infty}(I(x)+J_0(x)+\eta(x))=+\infty$ if $C$ is unbounded)
and so $(I+\eta)_{|\partial C}+\varphi$ is
weakly inf-compact in $\partial C$ since $\partial C$ is sequentially weakly closed (use also the Eberlein-Smulyan theorem). Therefore,
the assumptions of Theorem 2.3 are satisfied. Consequently, since $-S$ is convex and dense in $E^*$, there exists $\tilde\eta\in -S$, having the following property: 
for every function $J:C\to {\bf R}$ such that
$J_{|\partial C}-\varphi$ is constant in $\partial C$ and $I+J$ is lower semicontinuous and strictly convex, one has
$$\inf_{x\in C}(I(x)+J(x)+\tilde\eta(x))<\inf_{x\in \partial C}(I(x)+J(x)+\tilde\eta(x))\ .$$
In addition, let $J$ be G\^ateaux differentiable in $\hbox {\rm int}(C)$, with  $\lim_{\|x\|\to +\infty}{{I(x)+J(x)}\over {\|x\|}}=+\infty$ if $C$ is unbounded, and let $G:C\to {\bf R}$ be sequentially weakly lower semicontinuous in $C$
and G\^ateaux differentiable in $\hbox {\rm int}(C)$.
Fix $\sigma$ so that
$$\inf_{x\in C}(I(x)+J(x)+\tilde\eta(x))<\sigma<\inf_{x\in \partial C}(I(x)+J(x)+\tilde\eta(x))\ .\eqno{(2.2)}$$
 Since the set
$$\{x\in C : I(x)+J(x)+\tilde\eta(x)\leq\sigma\}$$
is sequentially weakly compact, in view
of Theorem 2.1 of [1], there exists $\epsilon>0$ such that, for every $\lambda\in [0,\epsilon]$,
 the restriction of the function $I+J+\tilde\eta+\lambda G$ to the set
$$\{x\in C : I(x)+J(x)+\tilde\eta(x)<\sigma\}$$
has a global minimum, say $\tilde x$. But, due to $(2.2)$, we have
$$\{x\in C : I(x)+J(x)+\tilde\eta(x)<\sigma\}\subseteq \hbox{\rm int}(C)$$
and, since $I+J+\tilde\varphi$ turns out to be continuous in $\hbox {\rm int}(C)$,  this implies that
$$I'(\tilde x)+J'(\tilde x)+\lambda G(\tilde x)+\tilde\eta=0\ ,$$
as claimed, with $\tilde\gamma=-\tilde\eta$.\hfill $\bigtriangleup$\par
\medskip
Now, we can give the\par
\medskip
{\it Proof of Theorem 1.1}. Let $S\subseteq E^*$ be a convex set dense in $E^*$ and let $\varphi:\partial X\to {\bf R}$.
Apply Theorem 2.5 with $I=0$. Let  $\tilde\gamma\in S$ be as in the conclusion of Theorem 2.5. Fix any lower semicontinuous and strictly convex function $J:C\to {\bf R}$, G\^ateaux differentiable in $\hbox {\rm int}(C)$, such that
$J_{|\partial C}-\varphi$ is constant in $\partial C$, with $\lim_{\|x\|_E\to +\infty}{{J(x)}\over {\|x\|}}=+\infty$ if $C$ is
unbounded. Now, fix any $G\in E^*$. Then, there exists $\epsilon>0$ such that the equation
$$J'(x)-\lambda G=\tilde\gamma$$
has a solution in $\hbox {\rm int}(C)$ for all $\lambda\in [0,\epsilon]$, and this means exactly that the set $\tilde\gamma$ is an algebraically
interior point of $J'(\hbox {\rm int}(C))$ \hfill $\bigtriangleup$\par
\medskip
In a finite-dimensional setting, another consequence of Theorem 2.5 is as follows:\par
\medskip
THEOREM 2.6. - {\it Let $E$ be a finite-dimensional real Banach space and let $C\subset E$ be a compact convex set with non-empty interior. \par
Then, for every function $\varphi:\partial C\to {\bf R}$, there exists $\tilde\gamma\in E^*$ having the following
property: for every lower semicontinuous stricly convex function
 $J:C\to {\bf R}$, G\^ateaux differentiable in $\hbox {\rm int}(C)$, such that $J_{|\partial C}-\varphi$
is constant in $\partial C$, and for every lower semicontinuous, bounded below and G\^ateaux differentiable function 
$H:\hbox {\rm int}(C)\to {\bf R}$,
 there exists $\epsilon>0$ such that, for each $\lambda\in [0,\epsilon]$, the equation
$$J'(x)+\lambda H'(x)=\tilde\gamma$$
has at least one solution in $\hbox {\rm int}(C)$.}
\smallskip
PROOF. Of course, since $C$ is compact, $\partial C$ is not convex. Fix $\varphi:\partial C\to {\bf R}$ and apply Theorem 2.5, with $I=0$. Let
$\tilde\gamma\in E^*$ be as in the conclusion of Theorem 2.5. 
Let $J:C\to {\bf R}$ be a lower semicontinuous strictly convex function, G\^ateaux
differentiable in $\hbox {\rm int}(C)$, such that $J_{|\partial C}-\varphi$
is constant in $\partial C$, and let $H:\hbox {\rm int}(C)\to {\bf R}$ be lower semicontinuous, bounded below and 
G\^ateaux differentiable. Now, consider the function $G:C\to {\bf R}$ defined by
$$G(x)=\cases{H(x) & if $x\in \hbox {\rm int}(C)$\cr & \cr
a & if $x\in \partial C$\ ,\cr}$$
where $a=\inf_CH$. Clearly, $G$ is lower semicontinuous in $C$. Consequently, there exists $\epsilon>0$ such that, for
each $\lambda\in [0,\epsilon]$, the equation
$$J'(x)+\lambda G'(x)=\tilde\gamma$$
has at least one solution in $\hbox {\rm int}(C)$, and the conclusion holds since $H=G$ in $\hbox {\rm int}(C)$.\hfill
$\bigtriangleup$
\medskip
Finally, we highlight the following consequence of Theorem 2.6:\par
\medskip
THEOREM 2.7. - {\it  Let $E$ be a finite-dimensional real Hilbert space and let $C\subset E$ be a closed ball, centered at $0$. \par
Then, for every function $\varphi:\partial C\to {\bf R}$, there exists $\tilde\gamma\in E^*$, having the following property: 
for every $P\in C^1(C)$ such that
$P'$ is Lipschitzian with Lipschitz constant $L\geq 0$, for every $\mu>L$, for every lower semicontinuous convex function 
$Q:C\to {\bf R}$, G\^ateaux differentiable in $\hbox {\rm int}(C)$, such that $(P+Q)_{|\partial C}-\varphi$ is constant in
$\partial C$, and for every lower semicontinuous, bounded below and  G\^ateaux differentiable function $H:\hbox {\rm int}(C)
\to {\bf R}$, there exists $\epsilon>0$ such that, for each $\lambda\in [0,\epsilon]$, the equation
$$\mu x+P'(x)+Q'(x)+\lambda H'(x)=\tilde\gamma$$
has at least one solution in $\hbox {\rm int}(C)$.}\par
\smallskip
PROOF. Fix $\varphi:\partial C\to {\bf R}$. Let $\tilde\gamma$ be as in the conclusion of Theorem 2.6.
Observe that, since $\mu>L$, the function ${{\mu}\over {2}}\|\cdot\|^2+P(\cdot)$ turns out to be
strictly convex since $E$ is a Hilbert space. Consequently, the function 
$$J(\cdot):={{\mu}\over {2}}\|\cdot\|^2+P(\cdot)+Q(\cdot)$$
is lower semicontinuous, strictly convex, G\^ateaux differentiable in $\hbox {\rm int}(C)$ and $J_{|\partial C}-\varphi$ is constant in
$\partial C$. Hence, the conclusion of Theorem 2.6 applies with such a function $J$ and we are done.\hfill $\bigtriangleup$\par
\bigskip
\bigskip
\bigskip
\bigskip
{\bf Acknowledgements.} The author has been supported by the Gruppo Nazionale per l'Analisi Matematica, la Probabilit\`a e 
le loro Applicazioni (GNAMPA) of the Istituto Nazionale di Alta Matematica (INdAM) and by the Universit\`a degli Studi di Catania, PIACERI 2020-2022, Linea di intervento 2, Progetto ”MAFANE”. The author would like to thank the referee for his/her comments and remarks.\par
\vfill\eject
\centerline   {\bf References}\par
\bigskip
\bigskip
\noindent                          
[1]\hskip 5pt B. RICCERI, {\it Nonlinear eigenvalue problems},  
in ``Handbook of Nonconvex Analysis and Applications'' 
D. Y. Gao and D. Motreanu eds., 543-595, International Press, 2010.\par
\smallskip
\noindent
[2]\hskip 5pt B. RICCERI, {\it On a minimax theorem: an improvement, a new proof and an overview of its applications},
Minimax Theory Appl., {\bf 2} (2017), 99-152.\par
\smallskip
\noindent
[3]\hskip 5pt B. RICCERI, {\it Multiplicity theorems involving functions with non-convex range}, Stud. Univ. Babe\c{s}-Bolyai Math., {\bf 68} (2023), 125-137.\par

\bigskip
\bigskip
Department of Mathematics and Informatics\par
University of Catania\par
Viale A. Doria 6\par
95125 Catania, Italy\par
{\it e-mail address}: ricceri@dmi.unict.it

\bye